\documentclass[11pt,a4paper]{article}
\usepackage[utf8]{inputenc}
\usepackage[T1]{fontenc}
\usepackage{amssymb}
\usepackage{amsthm}
\usepackage{color}
\definecolor{DarkOlive}{rgb}{0.1047,0.2412,0.0064}
\definecolor{FireBrick}{rgb}{0.5812,0.0074,0.0083}
\definecolor{RoyalBlue}{rgb}{0.0236,0.0894,0.6179}
\definecolor{RoyalGreen}{rgb}{0.0236,0.6179,0.0894}
\definecolor{RoyalRed}{rgb}{0.6179,0.0236,0.0894}
\definecolor{LightBlue}{rgb}{0.8544,0.9511,1.0000}
\definecolor{Black}{rgb}{0.0,0.0,0.0}
\definecolor{MidnightBlue}{rgb}{0.0035,0.0020,0.1363}
\definecolor{FireBrick3}{rgb}{0.5812,0.0074,0.0083}
\definecolor{FireBrick4}{rgb}{0.2156,0.0023,0.0035}
\definecolor{Blue2}{rgb}{0.0000,0.0000,0.8644}
\definecolor{Navy}{rgb}{0.0000,0.0000,0.1927}
\definecolor{MediumBlue}{rgb}{0.0000,0.0000,0.6179}
\usepackage[
        bookmarks=false,pdftitle={Representations of the same degree},
        colorlinks=true,backref=false,breaklinks=true,linkcolor=MediumBlue,
        citecolor=FireBrick3,filecolor=RoyalRed,
        urlcolor=Blue2]{hyperref}
\theoremstyle{plain}
\newtheorem{Thm}{Theorem}

\newtheorem{Pro}[Thm]{Proposition}

\theoremstyle{definition}

\newcommand{\mod}{ \textrm{ mod } }

\begin{document}
\title{Representations with the same degree}
\author{Frank Lübeck\thanks{This is a contribution to 
Project-ID 286237555 -- TRR 195 -- by the
German Research Foundation}}
\maketitle
\begin{abstract}
In this short note we show that every connected reductive simply-connected
algebraic group of rank $>1$ over the complex numbers has infinitely 
many pairs of 
irreducible representations which are not related by an automorphism of 
the algebraic group and which have the same degree.

This answers a question I was asked by J.~P.~Serre.
\end{abstract}

\section{Introduction and notation}\label{secIntro}

Let $G$ be a connected reductive algebraic group defined over the complex
numbers, and assume that $G$ is of simply-connected type. 
We will prove the following theorem.

\begin{Thm}
When the rank of $G$ is $>1$ then
$G$ has infinitely many pairs of finite dimensional irreducible
rational representations $\rho_1$ and $\rho_2$ such that there is no
automorphism $\alpha$ of $G$ with $\rho_2 = \rho_1 \circ \alpha$ and such
that $\rho_1$ and $\rho_2$ have the same degree.
\end{Thm}

We recall some facts and notation for which we refer to books on linear
algebraic groups, for example~\cite{spr,hum,milne}.

Let $B \subseteq G$ be a Borel subgroup and $T \subseteq B$ be a maximal torus
of $G$. To this pair one can associate a root datum $(X,\Phi,Y,\Phi^\vee)$.
If $l$ is the rank of $G$ then $X \cong \mathbb{Z}^l \cong Y$, and the
root datum comes with a non-degenerate pairing $\langle \cdot, \cdot \rangle
: X \times Y \to \mathbb{Z}$, $\Phi$ is a finite subset of $X$, called the
roots of $G$,  $\Phi^\vee$ is a finite subset of $Y$, called the coroots of
$G$, and there is a bijection $\Phi \to \Phi^\vee$, $\alpha \mapsto
\alpha^\vee$. The root datum determines $G$ up to isomorphism. 

Each root system is the disjoint union of irreducible root systems, and the
irreducible root systems are classified into four infinite series
$A_l$, $B_l$, $C_l$, $D_l$ and the exceptional cases $G_2$, $F_4$, $E_l$
with $6\leq l \leq 8$; the subscript is the rank $l$.

A root system contains a subset $\{\alpha_1, \ldots, \alpha_l\}$ of $l$
simple roots, the corresponding coroots $\{\alpha_1^\vee, \ldots, 
\alpha_l^\vee\}$ are a set of simple coroots. 

Our assumption that $G$ is simply-connected means that the simple coroots 
are a $\mathbb{Z}$-basis of the lattice $Y$. The elements
$\{\omega_1,\ldots, \omega_l\}$ of the dual basis of $X$ (with
respect to the pairing $\langle \cdot, \cdot \rangle$) are called the
fundamental weights. The positive coroots ${\Phi^\vee_+}$ are the 
elements $\alpha^\vee \in \Phi^\vee$ of the form $\alpha^\vee = 
b_1 \alpha_1^\vee + 
\ldots b_l \alpha_l^\vee$ with all $b_i \in \mathbb{Z}_{\geq 0}$
(then $\Phi^\vee = \Phi^\vee_+ \cup -\Phi^\vee_+$). So, for any weight
$\lambda = a_1 \omega_1 + \ldots + a_l \omega_l \in X$ we can compute
\[ \langle \lambda, \alpha^\vee\rangle = a_1 b_1 + \ldots a_l b_l. \]

The group $G$ is a direct product of simple simply-connected groups (which
correspond to the irreducible subsets of $\Phi$).

The irreducible rational representations  of $G$ are parameterized by the set
of dominant weights $\{\lambda = a_1 \omega_1 + \ldots + a_l \omega_l \mid\;
a_i \in \mathbb{Z}_{\geq 0}\}$.  We write $V(\lambda)$ for the module
corresponding to the dominant weight $\lambda$. We write $\rho = \omega_1 +
\ldots + \omega_l$ for the sum of the fundamental weights. 
The Weyl dimension formula 
descibes the degrees of the irreducible representations:

\[ \dim(V(\lambda)) = \frac{\prod_{\alpha^\vee \in \Phi^\vee_+} \langle
\lambda+\rho, \alpha^\vee\rangle}
{\prod_{\alpha^\vee \in \Phi^\vee_+} \langle
\rho, \alpha^\vee\rangle}. \] 

So, the only information we need for computing the degrees of the irreducible
representations of $G$ is the list of positive coroots, expressed as linear
combinations of the simple coroots.

Let $\sigma$ be an automorphism of the algebraic group $G$. Since the pair
$T \subseteq B$ is unique up to conjugation in $G$ there is an inner
automorphism $\tau$ of $G$ such that $\sigma \circ \tau$ fixes $T$ and $B$.
This automorphism then induces actions on $X$ and $Y$ which permute the sets
of simple coroots and of fundamental weights. It follows that $V(\lambda)$
and $V(\mu)$ cannot be related by an automorphism when the dominant weights
$\lambda$ and $\mu$ written as linear combinations of the fundamental weights
have different sets of coordinates. (See~\cite[Ch.11]{maltes} for a more
detailed description of the automorphisms of $G$.)

\section{Proof of Theorem~1}\label{secproof}

Let $G$ be a group of rank $l$ as in the previous section.

If for a dominant weight $\lambda$ we have $\dim(V(\lambda)) = d$ then it
is immediate from the Weyl dimension formula that $\dim(V(k (\lambda+\rho) -
\rho)) = d \cdot k^N$ for any $k \in \mathbb{Z}_{>0}$, where $N$ is the
number of positive (co-)roots of $G$.

Therefore, to show Theorem~1 it is sufficient to find for each $G$ of rank 
$l > 1$ one pair of dominant weights $\lambda$, $\mu$ which have
different sets of coordinates as linear combinations of fundamental weights
such that $\dim(V(\lambda)) = \dim(V(\mu))$.

To prove Theorem~1 we will describe an explicit pair  $\lambda$, $\mu$  in
each case.

\subsection{The case of rank $l=1$}

In this case $G \cong \mathop{SL}_2(\mathbb{C})$. Its coroot system of type
$A_1$ has only one simple
coroot and this is the only positive coroot. The Weyl dimension formula
yields $\dim(V(a \omega_1)) = (a+1)$ for all $a \in \mathbb{Z}_{\geq 0}$.
So, in this case all irreducible representations are determined by their
degree. The condition $l>1$ in Theorem~1  cannot be avoided.

From now we assume $l>1$.

\subsection{The case that all simple direct factors have rank $1$}

Then $G$ is a direct product of two or more copies of
$\mathop{SL}_2(\mathbb{C})$. Let $\omega_1, \ldots, \omega_l$ be the
fundamental weights (one for each simple factor).

The discussion in the previous subsection shows $V(5\omega_1)$ and
$V(\omega_1 + 2\omega_2)$ both have dimension $6$. 
This proves Theorem~1 in this case.

If $G$ has a simple direct factor of rank $l>1$ it is sufficient to show
Theorem~1 for this factor (extend by the zero weight corresponding to the 
the trivial representation on the other factors).

From now we assume that $G$ is a simple group of rank $l > 1$.

\subsection{The case of exceptional types, $A_2$ and $B_2$}

In these cases we just used the Weyl dimension formula to compute explicitly
the degrees of $V(\lambda)$ for a range of weights $\lambda$ with small
coordinates until we found examples that prove Theorem~1 in these cases.

Detailed descriptions of irreducible root systems, including a description
of all roots as linear combinations of simple roots, can be found 
in the tables~\cite[Planches I-IX]{bou}. In the results below we will 
use the ordering
of the simple coroots and corresponding fundamental weights $\omega_1,
\ldots, \omega_l$ as given in these tables.

\begin{sloppypar} 
For the computations we used our computer algebra package CHEVIE~\cite{chev}
which contains functions that compute these data.
\end{sloppypar}

Here are the smallest dimensional examples for the considered types (we have
checked all weights which yield smaller than the given dimensions).

\begin{Pro}
The following table gives for $G$ with irreducible root system of exceptional
type and of type $A_2$ and $B_2$ two dominant weights 
$\lambda$ and $\mu$ such that $d = \dim(V(\lambda)) = \dim(V(\mu))$.

\mbox{}

\quad\quad
  \begin{tabular}{l|l|l|l}
      \textbf{type} & \textbf{$\lambda$} & \textbf{$\mu$} & $d$ \\
      \hline
      $A_2$ & $\omega_1+2\omega_2$ & $4\omega_2$ & $15$ \\
      $B_2$ & $\omega_1+2\omega_2$ & $4\omega_2$ & $35$ \\
%% note that CHEVIE and Bourbaki use reverse ordering in type G_2
      $G_2$ & $3 \omega_1$ & $2 \omega_2$ & $77$\\
      $F_4$ & $\omega_1 + \omega_4$ & $2\omega_1$ & $1053$ \\
      $E_6$ & $2\omega_1$ & $\omega_3$ & $351$\\
      $E_7$ & $\omega_4+\omega_5$ & $2\omega_6+3\omega_7$ & $1903725824$\\
      $E_8$ & $\omega_1+\omega_3$ & $\omega_1+\omega_7+\omega_8$ & $8634368000$\\
  \end{tabular}

In case $E_6$ there is a  non-trivial graph automorphism of $G$ which shows
that $V(2\omega_6)$ and $V(\omega_5)$ also have degree $351$.
\end{Pro}

\section{Cases $A_l$, $B_l$, $C_l$, $D_l$ for large rank $l$}

As in the previous section we have computed many degrees of representations
in these cases up to rank $12$ and higher to look for patterns of 
dominant weights
that lead to the same degrees. Interestingly, we could detect
in all cases such patterns among  weights of the form $a\omega_1 + b\omega_2$.

In cases $A_l$, $B_l$, $D_l$ we even found examples with $a \in \{0,1\}$ for
all $l$, and the proof for general $l$ is straight forward.

The more interesting case are groups of type $C_l$ where the
smallest degrees which occur twice and the coefficients of the corresponding 
weights can be small for some $l$ and very large for other $l$.  
Here we reduce the existence of examples to the solvability of a Diophantine
equation, and show that solutions exist for all $l$.

\begin{Thm}
\begin{itemize}
\item[(a)] Let $G$ be of type $A_l$, $l\geq 3$. Then $V((l-1)\omega_2)$ and
$V(\omega_1 + (l-2)\omega_2)$ have the same degree
\[ \frac{(2l-1) \prod_{k=l+1}^{2l-2} k^2}{(l-1)!^2}. \]
(This is also correct for $l=2$ but then the two weights are related by a
graph isomorphism.)

\item[(b)] Let $G$ be of type $B_l$, $l\geq 3$. Then $V((2l-2)\omega_2)$ and
$V(\omega_1 + (2l-3)\omega_2)$ have the same degree
\[ \frac{3\cdot (4l-5)(6l-5)(6l-7) \prod_{k=2l}^{4l-6}k^2}{(2l-3)!^2}. \]

\item[(c)] Let $G$ be of type $D_l$, $l\geq 4$. Then $V((2l-3)\omega_2)$ and
$V(\omega_1 + (2l-4)\omega_2)$ have the same degree
\[ \frac{3 \cdot (3l-4)(3l-5)(4l-7)\prod_{k=2l-1}^{4l-8} k^2}{(l-2)^2
(2l-5)!^2}. \]

\item[(d)] Let $G$ be of type $C_l$, $l\geq 3$. Then for $a \geq 3$, $b \in
\mathbb{Z}_{\geq 0}$ the representations $V(a\omega_1 + b\omega_2)$ and
$V((a-2)\omega_1 + (b+1)\omega_2)$ have the same degree
\[ \frac{(a+1)(a+2b+2l-1) \prod_{k=a+b+2}^{a+b+2l-2} k \prod_{k=b+1}^{b+2l-3}
k}{(2l-1)!(2l-3)!} \]
if and only if there exists a positive integer $c$ such that 
\[ c^2 - (4l-5) a^2 = (2l-3)^2 \] 
and in that case $b = \frac{1}{2}(c+1-a-2l)$. 
For every $l \geq 3$ there exist infinitely many
pairs $(a,b)$ with this property.
\end{itemize}
\end{Thm}
%\begin{Prf}

\noindent\textbf{Proof. }\\[-1.5em]
\begin{itemize}
\item[(a)] We use the description in~\cite[Planche I]{bou}(II) which
describes all positive coroots as linear combinations of the simple ones:
\[ \beta_{i,j} = \sum_{k=i}^j\alpha_k^\vee \textrm{ for } 1 \leq i \leq j
\leq l.\]
Let $\lambda$ be a dominant weight of form $\lambda = a \omega_1 + b
\omega_2$. Evaluating the numerator in the Weyl dimension formula the only
factors which depend on $a$ and $b$ are:
\[ \langle \lambda+\rho, \beta_{1,1} \rangle = a+1,\]
\[ \langle \lambda+\rho, \beta_{1,j} \rangle = (a+1)+(b+1) + (j-2) = a+b+j
\textrm{ for } 2 \leq j \leq l,\]
\[ \langle \lambda+\rho, \beta_{2,j} \rangle = (b+1) + (j-2) = b-1+j
\textrm{ for } 2 \leq j \leq l.\]
Substituting in these terms $a = b = 0$ we get the  corresponding factors in 
the denominator of the Weyl dimension formula. All other factors are the
same in the numerator and denominator and cancel out.

Now we specialize to $(a,b) = (0,l-1)$ and $(a,b) = (1, l-2)$. The factors
involving $\beta_{1,j}$ for $j\geq 2$ are the same. And all but the
first or last factor involving $\beta_{2,j}$ for $j \geq 2$ are the same.
We get
\[ \frac{\dim(V((l-1)\omega_2))}{\dim(V(\omega_1+(l-2)\omega_2))} = 
\frac{1 \cdot (2l-2)}{2 \cdot (l-1)} = 1. \]
And for the degree of $V((l-1)\omega_2)$ we get
\[ \dim(V((l-1)\omega_2)) = \frac{\prod_{k = l+1}^{2l-1} k 
\prod_{k = l}^{2l-2} k}{\prod_{k = 2}^{l} k \prod_{k = 1}^{l-1} k} =
 \frac{(2l-1) \prod_{k=l+1}^{2l-2} k^2}{(l-1)!^2}.
\]

\item[(b)] Now let $G$ be of type $B_l$ with $l \geq 3$. Note that the
bijection $\Phi \to \Phi^\vee$, $\alpha \mapsto \alpha^\vee$, maps long roots
to short coroots and vice versa. The coroot system is of type $C_l$, so we
need to use~\cite[Planche III]{bou} in this case. We have the following 
positive coroots.
\[ \beta{ij} = \sum_{k=i}^{j-1} \alpha_k^\vee + 2 \sum_{k=j}^{l-1}
\alpha_k^\vee +\alpha_l^\vee \textrm{ for } 1\leq i \leq j \leq l.\]
\[ \gamma_{ij} = \sum_{k=i}^j \alpha_k^\vee \textrm{ for } 1\leq i \leq j
\leq l-1.\]

Let $\lambda = a \omega_1 + b \omega_2$. 
Only the following factors in the Weyl dimension formula for $V(\lambda)$
depend on $a$ and $b$:

\[ \langle \lambda+\rho, \beta_{1,1} \rangle = 2a + 2b +2l - 1,\]
\[ \langle \lambda+\rho, \beta_{1,2} \rangle = a + 2b +2l - 2,\]
\[ \langle \lambda+\rho, \beta_{1,j} \rangle = a + b +2l - j 
\textrm{ for } 3 \leq j \leq l,\]
\[ \langle \lambda+\rho, \beta_{2,2} \rangle = 2b +2l - 3,\]
\[ \langle \lambda+\rho, \beta_{2,j} \rangle = b +2l -1 - j 
\textrm{ for } 3 \leq j \leq l,\]
\[ \langle \lambda+\rho, \gamma_{1,1} \rangle = a + 1,\]
\[ \langle \lambda+\rho, \gamma_{1,j} \rangle = a + b + j
\textrm{ for } 2 \leq j \leq l-1,\]
\[ \langle \lambda+\rho, \gamma_{2,j} \rangle = b - 1 + j
\textrm{ for } 2 \leq j \leq l-1.\]

Specializing to $(a,b) = (0, 2l-2)$ and $(a,b)= (1, 2l-3)$ we see again that
most of the listed factors are the same in both cases. We get the following
quotient of the corresponding degrees:
\[ \frac{\dim(V((2l-2)\omega_2))}{\dim(V(\omega_1+(2l-3)\omega_2))} =
\frac{(6l-6)(6l-7)(4l-6)(3l-4)}{2(6l-7)(6l-9)(3l-4)(2l-2)} = 1.
\]

And for the degree we find after some simplifications
\[\dim(V((2l-2)\omega_2)) = 
\frac{3\cdot (4l-5)(6l-5)(6l-7) \prod_{k=2l}^{4l-6}k^2}{(2l-3)!^2}.
\]

\item[(c)] Now let $G$ be of type $D_l$. Here we need~\cite[Planche
IV]{bou}, but note that in this case the positive coroots are not
correctly described as linear combinations of simple ones in the reference. 
We correct this here. The positive coroots are:
\[ \beta_{ij} = \sum_{k=i}^{j} \alpha_k^\vee \textrm{ for }
1 \leq i \leq j \leq l-1,\]
\[ \eta_{i} = \sum_{k=i}^{l-2} \alpha_k^\vee + \alpha_l^\vee \textrm{ for }
1 \leq i \leq l-1,\]
\[ \gamma_{ij} = \sum_{k=i}^{j-1} \alpha_k^\vee + 2 \sum_{k=j}^{l-2}
\alpha_k^\vee +\alpha_{l-1}^\vee + \alpha_l^\vee.\]
The rest of the argument is completely analogous to~(a) and~(b). We
skip the details.

\item[(d)] Finally, we consider the case when $G$ is of type $C_l$ with $l
\geq 3$. Here, we describe the positive coroots using~\cite[Planche II]{bou} 
(type $B_l$).
\[ \beta_{ij} = \sum_{k=i}^j \alpha_k^\vee \textrm{ for } 1\leq i \leq l-1,\]
\[ \eta_i = \sum_{k=i}^l \alpha_k^\vee \textrm{ for } 1 \leq i \leq l,\]
\[ \gamma_{ij} = \sum_{k=i}^{j-1} \alpha_k^\vee + 2 \sum_{k=j}^l
\alpha_k^\vee \textrm{ for } 1 \leq i < j \leq l. \]

For $\lambda = a \omega_1 + b \omega_2$ the following factors in the
Weyl dimension formula depend on $a$ or $b$.

\[ \langle \lambda+\rho, \beta_{1,1}\rangle = a + 1,\]
\[ \langle \lambda+\rho, \beta_{1,j}\rangle = a + b + j \textrm{ for } 2\leq
j \leq l-1,\]
\[ \langle \lambda+\rho, \beta_{2,j}\rangle = b - 1 + j \textrm{ for } 2\leq
j \leq l-1,\]
\[ \langle \lambda+\rho, \eta_{1}\rangle = a + b + l,\]
\[ \langle \lambda+\rho, \eta_{2}\rangle = b + l - 1,\]
\[ \langle \lambda+\rho, \gamma_{1,2}\rangle = a + 2b + 2l - 1,\]
\[ \langle \lambda+\rho, \gamma_{1,j}\rangle = a + b + 2l + 1 - j 
\textrm{ for } 3\leq j \leq l,\]
\[ \langle \lambda+\rho, \gamma_{2,j}\rangle = b + 2l - j 
\textrm{ for } 3\leq j \leq l.\]

Taking the product of these factors and substituting $a=b=0$ in this product
yields numerator and denominator of $\dim(V(\lambda))$ as stated in the
theorem.

Now assume $a \geq 3$. Substituting in the factors above $a$ by $a-2$ and
$b$ by $b+1$, and cancelling all common factors we get
\[ \frac{\dim(V(a\omega_1 + b\omega_2))}{\dim(V((a-2)\omega_1 
+ (b+1)\omega_2))} = \frac{(a+1)(b+1)(a+b+2l-2)}{(a-1)(a+b+1)(b+2l-2)}.\]

Equating numerator and denominator and rearranging terms we see that this
quotient is $1$ if and only if
\[ b^2 + (a+2l-1) b + \frac{1}{2}((-2l+3) a^2 + (2l-1) a +(4l-4)) = 0.\]

Solving for $b$ yields
\[ b = \frac{1}{2}\left(-2l - a +1 \pm \sqrt{(4l-5) a^2 + (2l-3)^2}\right).\]

So, when $b$ is an integer then $(4l-5)a^2 + (2l-3)^2$ is the square $c^2$
of a positive integer $c$. Conversely, if $(4l-5)a^2 + (2l-3)^2 = c^2$ is the
square of a positive integer $c$ then 
\[b = \frac{1}{2}(-2l - a + 1 + c)\]
is also a positive integer: When $a$ is odd then $c$ is even, and when $a$ is
even then $c$ is odd; and $b$ is positive since $l\geq 3$ and $a\geq 3$.

We have shown that $\dim(V(a\omega_1 + b\omega_2)) = \dim(V((a-2)\omega_1 
+ (b+1)\omega_2))$ if and only if the Diophantine equation
\[ \quad\quad\quad  c^2 - (4l-5) a^2 = (2l-3)^2 \quad \quad\quad (*) \]
has a solution $(c,a)$ with $a \geq 3$. And in that case $b$ can be computed
as above.

The equation~(*) is a \emph{generalized Pell equation}. Obviously, 
$(c,a) = (2l-3, 0)$ is a solution of~(*), but it does not fulfill the
condition $a \geq 3$. 

Now we use solutions $(x,y) \in \mathbb{Z}_{>0}^2$ of the \emph{Pell
equation}
\[ \quad\quad x^2 - d y^2 = 1 \textrm{ with } d = (4l-5). \quad\quad (**)\]
First, note that $(4l-5) \equiv 3 \mod 4$ and so $d$ is not a square in
$\mathbb{Z}$. In that case equation~(**) has an infinite number of solutions,
see~\cite[Thm.7.25,7.26]{niv}, and they can be 
found as follows: $\sqrt{d}$ has a
continued fraction expansion of some period $r$. If $r$ is even, let $m=r-1$
and if $r$ is odd, let $m = 2r-1$, and compute the $m$-th continued fraction
approximation $\frac{x_1}{y_1}$ of $\sqrt{d}$. Then $(x_1 + y_1 \sqrt{d})
(x_1 - y_1 \sqrt{d}) = x_1^2 - d y_1^2 =1$ and $(x_1,y_1)$ is called the
fundamental solution of~(**); all solutions $(x_k,y_k)$ can be found by
computing $(x_1 + y_1 \sqrt{d})^k = x_k + y_k \sqrt{d}$ for $k \in
\mathbb{Z}_{>0}$.

Finally, note that the $(c_k,a_k) = ((2l-3)x_k, (2l-3) y_k)$ are 
infinitely many solutions of~(*), and almost all $a_k \geq 3$.
(In general these are not all solutions of~(*).) 
\hfill $\Box$
\end{itemize}
%\end{Prf}

\subsubsection*{Remarks}

\textbf{(1)}
For case $C_l$ in Theorem~3 we remark that the period length of the 
continued fraction of $4l-5$ can be large. For example for $l=159$ it is
$48$ and the smallest $(a,b)$ such that $V(a\omega_1+b\omega_2)$ and
$V((a-2)\omega_1+(b+1)\omega_2)$ have the same degree, which we construct in
our proof, is
\[(a,b) =
(613975804336172576474505,
7404460209629201092363289),\]
and the corresponding degree has $15728$ decimal digits. A brute force
search for solutions of equation~(*) in the previous proof for $l=159$ 
yields other
solutions, $(a,b) = (87, 902)$ is the smallest; but this is difficult to
predict and for other $l$ a brute force search did not yield any solution.

\noindent\textbf{(2)}
Integer sequences: The sequence of degrees in Theorem~3(a) for type
$A_l$, $l \geq 2$, occurs in numerous mathematical contexts, see the
entry A000891 of the Online Encyclopedia of Integer Sequences
(OEIS)~\cite{OEIS}. 

The sequence of degrees that occur for multiple dominant weights for 
type $G_2$ can be found in~\cite[A339248]{OEIS}.

The corresponding sequence for type $E_8$ can be found
in~\cite[A181746]{OEIS}.

\noindent\textbf{(3)}
Are there always degrees that occur more than two times?

P. Deligne showed me a proof (using
rational points of elliptic curves) that in type $A_2$ there exists for any
positive integer $m$ a degree that occurs for at least $m$ different
dominant weights.

Computing many examples we found that in type $A_l$ one easily finds degrees
that occur with high multiplicity. In the other types this seems to occur
less often. We only found some examples where degrees appeared more than two
times (say, $3$ to $7$ times). In cases $G_2$ and $E_8$ we did not find
any degree occuring more than two times (more precisely, we computed that
all degrees $< 10^{22}$ for type $G_2$ and all degrees $< 10^{60}$
for type $E_8$ occur at most for two dominant weights).

\newcommand{\etalchar}[1]{$^{#1}$}

\end{document}